\documentclass{amsart}
\usepackage{amssymb}

\newcommand{\new}[1]{}

\def\<{\langle}
\def\>{\rangle}

\def\<{\langle}
\def\>{\rangle}

\newcommand{\be}{\begin{equation}}
\newcommand{\ee}{\end{equation}}

\newcommand{\rf}[1]{(\ref{#1})}

\newcommand{\XX}{{\mathbf X}}

\newcommand{\sR}{{\mathbb R}}
\newcommand{\sC}{{\mathbb C}}
\newcommand{\sN}{{\mathbb N}}

\newcommand{\calA}{{\mathcal A}}

\newcommand{\calF}{{\mathcal F}}

\newcommand{\E}{{ E}}

\newcommand{\Var}{{\rm Var}}
\newcommand{\eps}{\varepsilon}

      \newtheorem{theorem}{Theorem}[section]
       \newtheorem{proposition}[theorem]{Proposition}
       
       \newtheorem{lemma}[theorem]{Lemma}

\theoremstyle{remark}
       \newtheorem{remark}{Remark}[section]
\theoremstyle{definition}

\author{
W{\l}odzimierz  Bryc
}
\thanks{\noindent Research partially supported by NSF
grant \#INT-0332062, by the C.P. Taft Memorial Fund, and University of Cincinnati's
Summer Faculty Research Fellowship Program}
\address{
Department of Mathematics,
University of Cincinnati,
PO Box 210025,
Cincinnati, OH 45221--0025, USA}
\email{Wlodzimierz.Bryc@UC.edu}

\author{Jacek Weso{\l}owski}
\address{ Faculty of Mathematics and Information Science\\
Warsaw University of Technology\\ pl. Politechniki 1\\ 00-661
Warszawa, Poland}
\email{wesolo@mini.pw.edu.pl}

\date{February 17, 2004. Revised: May 12, 2004. Corrected: November 18, 2004}

\title[$q$-Meixner processes]
{Conditional moments of $q$-Meixner processes}

\begin{document}

\tolerance=2000

\begin{abstract} We show that stochastic processes with linear conditional expectations and
quadratic conditional variances are Markov, and their transition probabilities are related to  a three-parameter family of
orthogonal polynomials which generalize the Meixner polynomials.
Special cases of these  processes are known  to arise from the
non-commutative generalizations of the L\'evy processes.
\end{abstract}
\maketitle
\section{Introduction}
\subsection{Motivation}
  It has been known since the work of
Biane \cite{Biane98} that every non-commutative process with free
increments gives rise to a classical Markov process, whose
transition probabilities  "realize" the non-commutative free
convolution of the corresponding measures.
It is natural to ask how to recognize in classical probabilistic terms which
 Markov processes  might arise from this construction. Unfortunately, the
non-commutative freeness seems to be poorly reflected in the
corresponding classical Markov process, which makes it hard to
answer this question. A more general framework might be less
constraining and easier to handle.

Non-commutative
processes with free increments can be thought as a special case corresponding
 to the value $q=0$ of
the more general class of $q$-L\'evy processes
\cite{Anshelevich01a}, \cite{Anshelevich03b}. Markov processes are known to arise in
 this more
general setting
 in two important cases: Bo{\.z}ejko, K{\"u}mmerer, and
Speicher, give explicit Markov transition probabilities for the
$q$-Brownian motion, see \cite[Theorem 1.10]{BKS97}, and
Anshelevich \cite[Corollary A.1]{Anshelevich03} proves the
corresponding result for the $q$-Poisson process.
Other $q$-L\'evy processes are still not
well understood, so it is not known whether Markov processes arise in the general case; for indications that
Markov property may perhaps fail, see \cite{Anshelevich03c}.

This paper arose as an attempt to better understand the emergence
of related Markov processes from probabilistic assumptions. We
define our class of processes by assuming that the first two
conditional moments are  given respectively by the generic linear
and quadratic expressions. Such assumptions are familiar from
L\'evy's characterization of the Wiener process as a martingale
and a quadratic martingale with continuous trajectories. For more
general processes the assumption of continuity of trajectories
fails, so we replace it  by  conditioning  with respect not only to the past,  but also to the future. This approach
originated with Pluci\'nska \cite{Plucinska83} who proved that
processes with linear conditional expectations and constant conditional
variances are Gaussian. Subsequent papers covered discrete Gaussian
sequences \cite{Bryc85b}, $L_2$-differentiable processes
\cite{Szablowski89}, Poisson process \cite{Bryc87c}, Gamma process \cite{Wesolowski89b}.
Weso{\l}owski \cite{Wesolowski93} unified several partial
results, identifying the general quadratic conditional variance
problem which characterizes the five L\'evy processes of interest
in this note: Wiener, Poisson, Pascal, Gamma, and Meixner. Our main result, Theorem \ref{T1}, extends \cite[Theorem
2]{Wesolowski93} to the more general quadratic conditional variances.
Similar analysis of stationary sequences in \cite{Bryc98} yields
the classical versions of the non-commutative  $q$-Gaussian
processes of \cite{BKS97}. Further contributions to the stationary
case can be found in \cite{MS02}.

Stochastic processes with linear conditional expectations and
quadratic conditional variances turn out to depend on three numerical parameters
$-\infty<\theta<\infty, \tau\geq 0$, and $-1\leq q\leq 1$. They are Markov,  and  arise from the non-commutative constructions, at least for
those
values of the parameters when such constructions are known. To point out the connection with the
orthogonal polynomials from which they are derived,  we call them  $q$-Meixner  processes.

When $q=1$,  the $q$-Meixner processes  have
independent increments and we recover the five L\'evy processes from \cite[Theorem 2]{Wesolowski93}. For other values of parameter $q$,
 we encounter several processes that arose in non-commutative
probability.\new{Modified!}
{   If $\tau=\theta=0$, we get the
 classical version of the
$q$-Brownian motion \cite{BKS97}.
If $\tau=0,
\theta\ne 0$ the $q$-Meixner processes arise as the classical version
from 
the $q$-Poisson
process defined in  \cite[Def. 6.16]{Anshelevich01a}.
}
When $q=0$ the $q$-Meixner
processes are related to the class of free L\'evy processes
considered by Anshelevich \cite{Anshelevich01}.

 The reasons why these special cases of $q$-Meixner processes should arise from the Fock space constructions are not clear to us.
 It is
not known whether the generic $q$-Meixner process arises as a
 classical version of a
non-commutative process, but the situation must be more complex. The
 connection  with the  $q$-Levy processes on the $q$-Fock space
 as defined in \cite{Anshelevich01a} fails for the following reason.
 In Proposition \ref{LM2}  below we
establish a polynomial martingale property \rf{proj} for all
 $q$-Meixner processes. But from Anshelevich
\cite[Appendix A.2]{Anshelevich03c} we know that a generic $q$-Levy
 process does
 not have martingale polynomials; the exceptions  are
$q=0,q=1$, the $q$-Poisson process, and the $q$-Brownian motion, and
 these are precisely the cases that we already mentioned above.

\subsection{Assumptions}
Throughout this paper $(X_t)_{t\geq 0}$ is a separable
square-integrable stochastic process, normalized so that for all
$t,s\geq 0$
\begin{equation}\label{cov}
E(X_t)=0,\: E(X_tX_s)=\min\{t,s\}.
\end{equation}
We are interested in
the processes
with linear conditional expectations and quadratic conditional
variances. More specifically,  we assume the following.

For all $0\leq s<t<u$,
\begin{equation}
\label{LR} E(X_t|\calF_{\leq s}\vee \calF_{\geq u})=aX_s+bX_u,
\end{equation}
where $a=a(s,t,u),b=b(s,t,u)$ are the deterministic functions of $s,t,u$,
and $\calF_{\leq s}\vee \calF_{\geq u}$ denotes the $\sigma$-field
generated by $\{X_t: t\in[0,s]\cup[u,\infty)\}$.


For ease of reference, we list the following trivial consequences of \rf{LR}.
From the form of the covariance it follows that
\begin{equation}\label{a+b}
a=\frac{u-t}{u-s},\:
b=\frac{t-s}{u-s}.
\end{equation}
Notice that from \rf{LR} we have

$$
E(E(X_t|\calF_s)-X_s)^2=E(E(E(X_t|\calF_{\leq s}\vee \calF_{\geq u})|\calF_s)-X_s)^2$$
$$=b^2E(E(X
_u-X_s|\calF_s))^2 \le (t-s)^2/(u-s).
$$
Passing to the
limit  as  $u\to \infty$ we see that
\begin{equation}\label{LR-}
E(X_t|\calF_{\leq s})=X_s
\end{equation}
  for $0\leq s\leq t$.
Similarly, taking $s=0$ in \rf{LR} we get
\begin{equation}\label{LR+}
E(X_t|\calF_{\geq u})=\frac{t}{u}X_u.
\end{equation}
Processes which satisfy condition \rf{LR} are sometimes called
harnesses, see \cite{Hammersley}, \cite{Williams}. We  assume in
addition that the conditional variance of $X_t$ given $\calF_{\leq
s}\vee\calF_{\geq u}$ is given by a quadratic expression in $X_s$,
$X_u$. Recall that the conditional variance of $X$  with respect
to a $\sigma$-field $\calF$ is defined as
$$\Var(X|\calF)=E(X^2|\calF)-\left(E(X|\calF)\right)^2.$$
For later calculations,
it is convenient to express this assumption as follows.

For all $0\leq s< t<u$,
\begin{eqnarray}
\label{QV} E(X_t^2|\calF_{\leq s}\vee \calF_{\geq u})
=AX_s^2+BX_sX_u+CX_u^2+D+\alpha X_s+\beta X_u,
\end{eqnarray}
where
$A=A(s,t,u),B=B(s,t,u),C=C(s,t,u),D=(s,t,u),\alpha=\alpha(s,t,u),
\beta=\beta(s,t,u)$ are the deterministic functions of $s,t,u$.

Since $X_0=0$, the coefficients
$a, A,B,\alpha$ are undefined at $s=0$. In some formulas for definiteness we
 assign these values by continuity.

It turns out that under mild assumptions,
the  functions $A,B,C,D,\alpha,\beta$, are determined uniquely as
explicit functions of $s,t,u$, up to some
numerical constants.
 The next assumption specifies two of these constants by requesting
 that $\Var(X_t|\calF_{\leq s})=const$ for all $0\leq s\leq t$.
 We use \rf{LR-} to state this assumption in the following  more explicit form.

\begin{equation}\label{CV}
\E(X_t^2|\calF_{\leq s})=X_s^2+t-s.
\end{equation}

Notice that equations \rf{LR-} and \rf{CV} imply that $\{X_t:t\geq
0\}$ and  $\{X_t^2-t:t\geq 0\}$ are martingales with respect to
the natural filtration $\calF_{\leq t}$; these two martingale
conditions (and continuity of trajectories) are the usual assumptions in the
L\'evy theorem.
%

\section{Conditional variances}
It is interesting to note that under mild assumptions, assumption
\rf{QV} can be written explicitly, up to some numerical constants.
Two of these numerical constants appear already under
 one-sided conditioning.
\begin{proposition}\label{T0}
Let $(X_t)_{t\geq0}$ be a separable square integrable stochastic
process which satisfies conditions \rf{cov}, \rf{LR}, and such
that $1, X_t, X_t^2$ are linearly independent for all $t>0$.
 If for every
$0< t<u$ the conditional expectation $E(X_t^2|{\calF}_{\ge u})$
is a quadratic expression in variable $X_u$, then
there are constants $\tau\in[0,\infty]$ and $\theta\in{ \sR}$ such that
\begin{equation}\label{RV-**}
\Var(X_t|\calF_{\geq u})=\left\{\begin{array}{ll}
\frac{t(u-t)}{u+\tau}\left(\tau \frac{X_u^2}{u^2}+\theta \frac{X_u}{u}+1\right)&
\mbox{if $\tau<\infty$,} \\
t(u-t)\left(\frac{X_u^2}{u^2}+\theta \frac{X_u}{u}\right)& \mbox{if $\tau=\infty$.} \\
\end{array}\right.
\end{equation}
for all $0\le
t<u$.
\end{proposition}
\begin{proof} By assumption, for any $0<s<t$
\begin{equation}\label{gf}
E(X_s^2|\calF_{\ge t})=m(s,t)X_t^2+n(s,t)X_t+o(s,t)\;,
\end{equation}
where $m$, $n$, $o$ are some functions.

On the other hand from \rf{LR+} we get
$$
E(X_sX_t|{\calF}_{\ge u})=E(E(X_s|{\calF}_{\ge t})X_t|{\calF}_{\ge u})=
\frac{s}{t}E(X_t^2|{\calF}_{\ge u}),
$$
and from \rf{LR} we get
$$
E(X_sX_t|{\calF}_{\ge u})=E(X_sE(X_t|\calF_{\leq s}\vee \calF_{\geq u})|{\calF}_{\ge u})
$$
$$=\frac{u-t}{u-s}E(X_s^2|{\calF}_{\ge
u})+\frac{t-s}{u-s}X_uE(X_s|{\calF}_{\ge u})
$$
$$=\frac{u-t}{u-s}E(X_s^2|{\calF}_{\ge
u})+\frac{(t-s)s}{(u-s)u}X_u^2\;.
$$
Combining the above two formulas we have
\begin{equation}\label{rel}
\frac{s}{t}E(X_t^2|{\calF}_{\ge u})=
\frac{u-t}{u-s}E(X_s^2|{\calF}_{\ge u})+\frac{(t-s)s}{(u-s)u}X_u^2.
\end{equation}
Now we substitute the conditional moments from (\ref{gf}) into   (\ref{rel}),  getting
$$
\frac{s}{t}\left(m(t,u)X_u^2+n(t,u)X_u+o(t,u)\right)$$
$$=
\frac{u-t}{u-s}\left(m(s,u)X_u^2+n(s,u)X_u+o(s,u)\right)+\frac{(t-s)s}{(u-s)u}X_u^2\;.
$$
Recall that $1,X_u,X_u^2$ are linearly independent.
Comparing the coefficients of respective powers of $X_u$ we obtain
$$
\frac{s}{t}m(t,u)=\frac{u-t}{u-s}m(s,u)+\frac{(t-s)s}{(u-s)u}\;,
$$
$$
\frac{s}{t}n(t,u)=\frac{u-t}{u-s}n(s,u)\;,\;\;\;\frac{s}{t}o(t,u)=\frac{u-t}{u-s}o(s,u)\;.
$$
The first equation leads to
$$
\left(\frac{m(t,u)}{t}-\frac{1}{u}\right)\frac{1}{u-t}=
\left(\frac{m(s,u)}{s}-\frac{1}{u}\right)\frac{1}{u-s},
$$
and hence
$$
m(t,u)=\frac{t}{u}+t(u-t)i(u)
$$
for some function $i:\sR\to\sR$.
The next two equations give
$$
n(t,u)=t(u-t)j(u)\;\;\;\mbox{ and }\;\;\;o(t,u)=t(u-t)k(u)\;
$$
for some functions $j,k:\sR\to\sR$.
Thus from \rf{gf} we get
$$
E(X_t^2|{\calF}_{\ge
u})=\left(\frac{t}{u}+t(u-t)i(u)\right)X_u^2+t(u-t)j(u)X_u+t(u-t)k(u)\;.
$$
Taking the expectations of both sides we get
$
t=t+tu(u-t)i(u)+t(u-t)k(u),
$ so $k(u)=-ui(u)$. Finally we have
\begin{equation}\label{rep1}
E(X_t^2|{\calF}_{\ge
u})=\frac{t}{u}X_u^2+t(u-t)\left[i(u)(X_u^2-u)+j(u)X_u\right]\;.
\end{equation}

To identify the functions $i$ and $j$ we fix  $s<t<u$ and insert
(\ref{rep1}) into the formula
$$
E(X_s^2|{\calF}_{\ge u})=E(E(X_s^2|{\calF}_{\ge t})|{\calF}_{\ge u})\;.
$$

This gives
$$
\frac{s}{u}X_u^2+s(u-s)\left[i(u)(X_u^2-u)+j(u)X_u\right]$$
$$=
E\left(\left.\frac{s}{t}X_t^2+s(t-s)
\left[i(t)(X_t^2-t)+j(t)X_t\right]\right|{\calF}_{\ge u}\right)$$ $$=
\frac{s}{t}\left\{\frac{t}{u}X_u^2+t(u-t)\left[i(u)(X_u^2-u)+j(u)X_u\right]\right\}$$$$+
s(t-s)i(t)\left\{\frac{t}{u}X_u^2+t(u-t)\left[i(u)(X_u^2-u)+j(u)X_u\right]\right\}$$
$$+s(t-s)j(t)\frac{t}{u}X_u-st(t-s)i(t)\;.
$$
Comparing the coefficients of respective powers of $X_u$ we obtain
\begin{equation}\label{i}
ui(u)=ti(t)+(u-t)ti(t)ui(u)\;,
\end{equation}
\begin{equation}\label{j}
uj(u)=tj(t)+(u-t)ti(t)uj(u)\;.
\end{equation}
If $i$ is non-zero for all $t>0$ then (\ref{i}) gives
$
\frac{1}{ti(t)}+t=\frac{1}{ui(u)}+u
$.  This means that
$\frac{1}{ti(t)}+t=-\tau$ for some constant $\tau$,
and  $\tau\geq 0$ since $1/i(t)$  cannot vanish for any $t>0$.
Hence
\begin{equation}\label{i+} i(t)=-\frac{1}{t(t+\tau)}\;.
\end{equation}
Using this in \rf{j} we get
$u(u+\tau)j(u)=t(t+\tau)j(t)$. Thus
$$j(t)=\frac{\theta}{t(t+\tau)}$$
for some real constant $\theta$. We get
\begin{equation}\label{fut2}
E(X_s^2|{\calF}_{\ge
t})=\frac{s(s+\tau)}{t(t+\tau)}X_t^2+\frac{s(t-s)}{t(t+\tau)}\theta
X_t+\frac{s(t-s)}{t+\tau}\;.
\end{equation}
Suppose now that $i(t)=0$ for some $t>0$. Then \rf{i} implies
that $i$ is a zero function, corresponding to $\tau=\infty$ in \rf{i+}. In this case \rf{j} leads to
$uj(u)=tj(t)$, which means that $j(t)=\theta/t$ for some real
number $\theta$. Thus in this case
\begin{equation}\label{fut1}
E(X_s^2|{\calF}_{\ge t})=\frac{s}{t}X_t^2+\frac{s(t-s)}{t}\theta
X_t\;.
\end{equation}

\end{proof}

Notice that taking the expected value of both sides of \rf{QV}, we get
 a trivial relation
\begin{equation}
\label{0} t-A s-Cu=B s+D,
\end{equation}
valid for all $0\leq s<t<u$.
We need additional relations between the coefficients in \rf{QV}.

\begin{lemma}\label{L000}
Let $(X_t)_{t\geq0}$ be a separable square integrable stochastic
process which satisfies conditions \rf{cov}, \rf{LR}, and
such that $1, X_t, X_t^2$ are linearly independent for all $t>0$.
Suppose that condition \rf{QV} holds with $D(s,t,u)\ne 0$ for  all
$0\leq s<t<u$. Then the conditional expectation $E(X_u^2|{\calF}_{\le t})$ is
quadratic in $X_t$ for any $0\le t<u$. Moreover,
\begin{equation}\label{E1}
E(X_u^2-u|\calF_{\leq s})=\left(1+\frac{A+B+C-1}{b-C}\right)(X_s^2-s)+\frac{\alpha+\beta}{b-C}X_s.
\end{equation}

\end{lemma}
\begin{proof}
Equation \rf{LR-} implies that $E(X_t^2 |\calF_{\leq s})=E(X_tX_u |\calF_{\leq s})$, so from \rf{LR} we get
 $$E(X_t^2 |\calF_{\leq s})=aX_s^2+bE(X_u^2 |\calF_{\leq s}).$$
From  \rf{QV} we get
$$E(X_t^2 |\calF_{\leq s})=AX_s^2+BX_s^2+CE(X_u^2|\calF_{\leq s})+(\alpha+\beta)X_s+D.$$
Notice that this implies $C\ne b$. Indeed, if $C=b$ then subtracting the equations we get a quadratic
equation for $X_s$. If this equation is non-trivial, then  $1,X_s,X_s^2$ are linearly dependent.
So  the coefficients in the quadratic  equation must all be
zero; in particular,  $D=0$, contradicting the assumption.

Since $C\ne b$, we can solve the equations for $E(X_t^2 |\calF_{\leq s})$ and
$E(X_u^2 |\calF_{\leq s})$. Using \rf{0}, we get \rf{E1} after a calculation.
\end{proof}

\begin{lemma}\label{L001}
Let $(X_t)_{t\geq0}$ be a separable square integrable stochastic
process which satisfies conditions \rf{cov}, \rf{LR}, \rf{CV} and
such that $1, X_t, X_t^2$ are linearly independent for all $t>0$.
Suppose that condition \rf{QV} holds with $D(s,t,u)\ne 0$ for  all
$0\leq s<t<u$. Then the conditional expectations
$E(X_s^2|{\calF}_{\ge t})$ 
are quadratic in $X_t$ for any $0\le s<t$. Moreover, there are
constants  $0\leq \tau<\infty,-\infty<\theta<\infty$
such that 
\rf{RV-**} holds true, and 
the parameters in \rf{QV}, evaluated at $0\leq s<t<u$, satisfy the following equations.
\begin{eqnarray}
\label{1}
A+B+C&=&1,\\
{As^2+Bsu+Cu^2-t^2}{}&=&\tau D \label{tau},\\
\label{alpha/D'}\label{3}
{s\alpha+u\beta}{}&=&\theta D,\\
\label{beta/D}
\alpha+\beta&=&0.
\end{eqnarray}
\end{lemma}
\begin{proof}
Comparing the coefficients in \rf{CV} and \rf{E1}, we  get \rf{1}, and \rf{beta/D}.

Setting $s=0$ in \rf{QV} we see that $E(X_t^2|{\calF}_{\ge u})$ is quadratic in $X_u$. Thus
Proposition \ref{T0} implies that \rf{RV-**} holds true. Notice that since $D(0,t,u)\ne 0$, we must have
$\tau<\infty$, so \rf{fut2} holds. We use the latter
in
$$E(X_t^2 |\calF_{\geq u})=
AE(X_s^2|\calF_{\geq u})+\frac{s}{u} BX_u^2+CX_u^2+(\frac{s}{u}\alpha+\beta)X_u+D,$$
which follows from \rf{QV}.
We get \rf{tau} from the comparison of the quadratic terms,
and  \rf{alpha/D'}
from the comparison of the linear terms.
\end{proof}
For future reference we state the following.
\begin{remark}\label{RRR}
The system of equations \rf{0}, \rf{1}, \rf{3},
 \rf{tau}, \rf{alpha/D'}, \rf{beta/D} has the solution
\begin{eqnarray}
 \alpha  &=& D \frac{-\theta}{u-s}, \label{a}
\\
\beta &=& D\frac{\theta}{ u-s}, \label{b}
\\A &=&\frac{ta}{s}- D\frac{u+\tau }
        {s(u-s)}, \label{A}
\\
B &=& D\frac{ s + u +
       \tau  }{s ( u - s )}-\frac{u-s}{s}ab, \label{B} \\
C &=& b-
    D\frac{ 1 }{u - s}. \label{C}
\end{eqnarray}
\end{remark}

We need the following  version of {\cite[Theorem 2]{Wesolowski93}}.
\begin{proposition}\label{L1}
 Let $(X_t)_{t\geq0}$ be a
separable square integrable stochastic process which satisfies
conditions \rf{cov}, \rf{LR}, \rf{QV}, and \rf{CV}. Suppose that
the coefficient $D$ in \rf{QV} satisfies $D(s,t,u)\ne 0$
for  all 
$0\leq s<t<u$,
and that $1, X_t, X_t^2$ are linearly independent for all $t>0$.
Then $E(|X_t|^p)<\infty$ for all $p\geq 0$.
%

Moreover, if $(X_t)$ and $(Y_t)$ satisfy these assumptions
with the same coefficients in \rf{QV}, then the joint moments of both processes are equal,
$$
E(X_{t_1}^{n_1}X_{t_2}^{n_2}\dots X_{t_k}^{n_k})=E(Y_{t_1}^{n_1}Y_{t_2}^{n_2}\dots Y_{t_k}^{n_k})
$$
for all $t_1,t_2,\dots,t_k>0, n_1,n_2,\dots,n_k\in\sN, k\in \sN$.
\end{proposition}
\begin{proof}
Fix $s<t$ and let $\{t_k:k\geq 0\}$ be an arbitrary infinite
strictly increasing sequence which contains $s$ and $t$ as consecutive elements, say
$s=t_{N}, t=t_{N+1}$ for some $N\in\sN$.

We apply \cite[Theorem 2]{Wesolowski93} to the sequence
$\xi_k=X_{t_k}$. Of course,
$\sigma(\xi_1,\dots,\xi_k,\xi_{k+1})\subset \calF_{\leq
t_{k-1}}\vee \calF_{\geq t_{k+1}}$. Therefore, conditions
\rf{LR-}, \rf{LR}, \rf{CV}, and \rf{QV} imply \cite[(6), (7), (8),
and (9)]{Wesolowski93}, respectively. Since $corr (\xi_{k-1},
\xi_k)=\sqrt{t_{k-1}/t_k}\ne 0,\pm1$,  the assumption
\cite[(10)]{Wesolowski93} holds true, too. Finally, notice that
Weso{\l}owski's $\alpha_k=1$, and his
$\underline{a}_k=C(t_{k-1},t_k,t_{k+1})\ne
a_k=b(t_{k-1},t_k,t_{k+1})$ because from \rf{C} we see that $D \ne
0$ if and only if $C\ne b$. Thus \cite[(11)]{Wesolowski93} hold
true. From
 \cite[Theorem 2]{Wesolowski93} we see that $E(|X_t|^p)<\infty$ for all $p>0$, and
that for $n=1,2\dots$, the conditional moment
$E(X_t^n|X_{t_{1}},\dots,X_{t_{N-1}},X_s)$
 is a unique polynomial of degree $n$ in the variable $X_s$.

If two processes satisfy the assumptions, then the conditional moments of both processes can be
expressed as polynomials with the same coefficients. This implies that all joint moments of the
processes are equal.
\end{proof}

Next we give the general form of the conditional variance under
the two-sided conditioning.

\begin{proposition}\label{iii'}   Let $(X_t)_{t\geq0}$ be a
separable square integrable stochastic process which satisfies
conditions \rf{cov}, \rf{LR}, \rf{QV}, and \rf{CV}. Suppose that
the coefficient $D$ in \rf{QV} satisfies $D(s,t,u)\ne 0$
for  all 
$0\leq s<t<u$, and that $1, X_t, X_t^2$ are linearly independent
for all $t>0$. Then there are parameters
$-\infty<\theta<\infty$, and $0\leq \tau<\infty$  such that the
first part of \rf{RV-**} holds true. In addition, there exists
$-1<q\leq 1$ such that
\begin{eqnarray}\label{q-Var}
&\Var(X_t|\calF_{\leq s}\vee\calF_{\geq u})=&\\
\nonumber & \displaystyle \frac{(u-t)(t-s)}{u+\tau-qs}\left(
(1-q)\frac{(X_u-X_s)(sX_u-uX_s)}{(u-s)^2}
+\tau\frac{(X_u-X_s)^2}{(u-s)^2}+\theta\frac{X_u-X_s}{u-s}+1
\right).&
\end{eqnarray}
\end{proposition}
\begin{proof}
By 
Proposition \ref{L1}, all moments of $X_t$ are finite.
Fix $s<t$. Then from \rf{LR+} and \rf{CV} we get
$
\frac{s}{t}E(X_t^3)= E(X_t^2X_s)= EX_s^3,
$
so $ E(X_t^3)/t$ does not depend on $t>0$.
On the other hand, from \rf{fut2} we get
$$E(X_s^3)=E(X_s^2X_t)=\frac{s(s+\tau)}{t(t+\tau)}E(X_t^3)+\frac{s(t-s)}{t+\tau}\theta.
$$
Hence
\begin{equation}\label{m3}
 EX_t^3=t \theta.
\end{equation}

Similarly, from \rf{CV} we get
$$
E(X_t^2X_s^2)= EX_s^4 +s (t-s),
$$
and from \rf{RV-**} we get
$$
E(X_t^2X_s^2)=
\frac{s(s+\tau)}{t(t+\tau)} EX_t^4+  \theta \frac{s(t-s)}{t(t+\tau)} E(X_t^3)+s t \frac{t-s}{t+\tau}.
$$
Using \rf{m3}
we get after a calculation  that
$
\frac{E(X_s^4)-s (s+\theta^2)}{s(s+\tau)}
$
does not depend on $s$. Thus
\begin{equation}\label{tmp1}
E(X_t^4)=(1+q)t(t+\tau)+t(t+\theta^2)
\end{equation}
for some constant $q\in\sR$.

A calculation gives
$$E(X_t-X_s)^2=t-s,\; E(X_t-X_s)^3=\theta(t - s),
$$
and
$$E(X_t-X_s)^4=
( t-s)\left(6s+\theta^2-\tau+(2+q)(t+\tau-3s)\right).$$

Since the determinant
$$
\frac{1}{(t-s)^2}\det
\left[ {\begin{array}{*{20}c}
   1 & {E(X_t-X_s)} & {E((X_t-X_s)^2)}  \\
   {E(X_t-X_s)} & {E((X_t-X_s)^2)} & {E((X_t-X_s)^3)}  \\
   {E((X_t-X_s)^2)} & E((X_t-X_s)^3) & {E((X_t-X_s)^4)}  \\
 \end{array} } \right]$$
$$=
 q\left( t + \tau-3s   \right) +s + t  + \tau $$
is non-negative, taking  $s=t-1$ and $t\to\infty$, we get $q\leq 1.$ Since
$1,X_t,X_t^2$ are linearly independent, the  determinant evaluated
at $s=0$ must be strictly positive, see \cite[pg. 19]{Chihara}.
This shows that $q> -1$.

It remains to determine the coefficients in \rf{QV}. Fix  $s<t<u$.
Comparing the two representations of $\E(X_tX_u^2|\calF_{\leq s})$
as
$$\E(\E(X_t|\calF_{\leq s}\vee \calF_{\geq u})X_u^2|\calF_{\leq s})=
\E(X_t \E(X_u^2|\calF_{\leq t})|\calF_{\leq s}),$$
and the similar two expressions for $\E(X_t^2X_u|\calF_{\leq s})$,
 we get two different expressions for $\E(X_t^3|\calF_{\leq s})$.
Equating them, we get
$$
a X_s^3+b \E(X_u^3|\calF_{\leq s})$$$$=
A X_s^3+B X_s^3+B X_s (u-s)+C \E(X_u^3|\calF_{\leq s})+D X_s+(\alpha+\beta) X_s^2+\beta(u-s).
$$
We can solve this  equation for $E(X_u^3|\calF_{\leq s})$, as
\rf{C} implies that $C\ne b$. Using \rf{b} and \rf{1}, the answer
simplifies to
$$
E(X_u^3|\calF_{\leq s})=X_s^3+\frac{B(u-s)+D}{b-C}X_s+(u-s)\frac{\beta}{b-C}.
$$
From this we get
$$\frac{s}{u}E(X_u^4)=E(X_sX_u^3)=E(X_s^4)+\frac{B(u-s)+D}{D}(u-s)s .$$
Substituting \rf{tmp1}
we deduce the following equation
\begin{equation}\label{D}
\frac{(u-s)B}{D}=1+q.
\end{equation}

Solving together equations \rf{1}, \rf{tau}, \rf{3},  \rf{beta/D}, and \rf{D} for
$A,B,C,D,\alpha,\beta$ we obtain \rf{q-Var}.
\end{proof}

\begin{remark}
Solving together equations \rf{1}, \rf{tau}, \rf{3},  \rf{beta/D}, and \rf{D} for
$A,B,C,D,\alpha,\beta$ we get
\begin{eqnarray}
\label{A'} A&=&\frac{u-t}{u-s} \times \frac{u+\tau-qt}{u+\tau-qs},\\
\label{B'} B&=&(1+q)\frac{t-s}{u-s}\times\frac{u-t}{u+\tau-qs},\\
\label{C'} C&=&\frac{t-s}{u-s}\times\frac{t+\tau-qs}{u+\tau-qs},\\
\label{D'} D&=&\frac{(u-t)(t-s)}{u+\tau-qs},\\
\label{a'} \alpha&=&-\theta \frac{(u-t)(t-s)}{(u-s)(u+\tau-qs)},\\
\label{b'} \beta&=&\theta \frac{(u-t)(t-s)}{(u-s)(u+\tau-qs)}.
\end{eqnarray}
\end{remark}
\new{New}{
\begin{remark} From the formula for $E(X_t-X_s)^4$ we see
that except for the case $q=1$, the increments of the process $X_t$ are not stationary.
For $\tau=0$, the increments of the corresponding non-commutative processes are stationary,
but this property is not inherited by the
classical version.
\end{remark}
}
\section{$q$-Meixner Markov processes}
 We use the standard notation
\begin{eqnarray*}
{[n]_{q}} &=&1+q+\dots +q^{n-1}, \\
{[n]_{q}!} &=&[1]_{q}[2]_{q}\dots [n]_{q},
\\ \left[ \begin{array}{c}n \\ k \end{array}\right] _{q} &=&\frac{[n]_{q}!}{[n-k]_{q}![k]_{q}!},
\end{eqnarray*}%
with the usual conventions $[0]_{q}=0,[0]_{q}!=1$.
For fixed real parameters $x, s, t, q,\theta,\tau$, define the polynomials  $Q_n$ in variable $y$
 by the three step recurrence
\begin{eqnarray}\label{Q-rec-G}
&\displaystyle yQ_n(y|x) = Q_{n+1}(y|x) & \\ \nonumber+ &(\theta
[n]_q+x q^n)Q_n(y|x)+(t-sq^{n-1}+\tau[n-1]_q)[n]_q Q_{n-1}(y|x),&
\end{eqnarray}
where $n\geq 1$, and $Q_{-1}(y|x)=0$, $Q_0(y|x)=1$, so $Q_1(y|x)=y-x$.
It is well known that  such polynomials are orthogonal with respect to
a probability measure if the last coefficient
 of the three step recurrence is 
positive, see \cite[Theorem I.4.4]{Chihara}.
Therefore,  \rf{Q-rec-G} defines a probability measure
whenever $x,\theta\in\sR, 0<s<t, \tau\geq 0, -1\leq q\leq 1$.
Moreover, in this case
\begin{equation}\label{pre-Carleman}
\sum_n \frac{1}{\sqrt{(t-sq^{n-1}+\tau[n-1]_q)[n]_q}}=\infty,
\end{equation}
so from Carleman's criterion (see \cite[page 59]{Shohat-Tamarkin}),
 this measure is unique. We denote this unique
probability measure by
$\mu_{x,s,t}(dy)$.

Of course,  $\mu_{x,s,t}(dy)=\mu_{x, s, t, q,\theta,\tau}(dy)$
depends on all the parameters of the recurrence \rf{Q-rec-G}.
It is worth noting explicitly that  if  $q=-1$ then  $[2]_q=0$, so
$\mu_{x,s,t}(dy)$ is supported on two points only.
In general, more explicit expressions for $\mu_{x,s,t}(dy)$ can perhaps be derived from
\cite[Theorem 2.5]{Askey-Wilson85} by taking their parameters $b=c=0$,
$ad=-(s(1-q)+\tau)/(t(1-q)+\tau)$, $a+d=((1-q) x-1)/\sqrt{t+\tau/(1-q)}$.

If we need to indicate the dependence of the polynomials in
\rf{Q-rec-G} on the additional parameters in the recurrence
\rf{Q-rec-G}, we  write  $Q_n(y|x,s,t)$.

We will need two algebraic identities;
the first one resembles \cite[(2.3)]{Al-Salam76} but is in fact different;
the second one is a slight generalization  of \cite[Theorem 1]{Bryc-Matysiak-Szablowski}.
\begin{lemma}\label{QQQ}
For every $x,y,z\in\sR$, $n\in\sN$,  and $0\leq s\leq t \leq u$ we \new{$0\leq s\leq t \leq u$ zbedne?}
have
\begin{eqnarray}\label{I}
\displaystyle  Q_n(z|x,s,u) 
= \displaystyle
\sum_{k=0}^n \left[ \begin{array}{c}n \\ k \end{array}\right]_q Q_{n-k}(y|x,s,t)Q_k(z|y,t,u).
\end{eqnarray}
\new{added $\left[ \begin{array}{c}n \\ k \end{array}\right]$ }
Furthermore, 
\begin{eqnarray}\label{BMS}
&\displaystyle Q_n(z|y,t,u)&
\\ \nonumber
=&\displaystyle \sum_{k=1}^n
\left[ \begin{array}{c}n \\ k \end{array}\right]_q Q_{n-k}(0|y,t,0)\left(Q_k(z|0,0,u)-Q_k(y|0,0,t)\right).&
\end{eqnarray}
\end{lemma}
\begin{proof}
Consider first the case $|q|<1$. It is easy to check
by $q$-differentiation with respect to $\zeta$  that
the generating function
$$\phi(\zeta,y,x,s,t)=\sum_{n=0}^\infty \frac{\zeta^n}{[n]_q!}Q_n(y|x,s,t)$$
of the polynomials $Q_n$  is given by
$$
\phi(\zeta,y,x,s,t)=\prod_{k=0}^\infty\frac{1+\theta \zeta  q^k-(1-q)x\zeta q^k+((1-q)s+\tau)\zeta^2 q^{2k}}
{1+\theta \zeta q^k-(1-q)y\zeta q^k+((1-q)t+\tau)\zeta^2 q^{2k}}.
$$
For details, see \cite{Al-Salam76}. Notice that for $|q|<1$, the series defining $\phi(\zeta,y,x,s,t)$ converges for all
 $|\zeta|$ small enough.
Indeed, from \rf{Q-rec-G} we get by induction $|Q_{n+1}|\leq C^n$ with
$C=\max\{1,(|x|+|y|+|\theta|+\tau+t+s)/(1-|q|)^2\}$.

Therefore,
\begin{equation}\label{product}
\phi(\zeta,z,x,s,u)=\phi(\zeta,y,x,s,t)\phi(\zeta,z,y,t,u),
\end{equation}
which implies \rf{I} for all $n\geq 0$ and $|q|<1$. Since \rf{I}
is an identity between the polynomial expressions in variables
$z,y,q$, it must hold for all  $q$.

Since ${1}/{\phi(\zeta,y,x,s,t)}=\phi(\zeta,x,y,t,s)$, from
\rf{product}  we get
$$
\phi(\zeta,z,y,t,u)=\frac{\phi(\zeta,z,x,s,u)}{\phi(\zeta,y,x,s,t)}$$
$$
=1+\frac{1}{\phi(\zeta,y,x,s,t)}\left(\phi(\zeta,z,x,s,u)-\phi(\zeta,y,x,s,t)\right)
$$
$$
=1+\phi(\zeta,x,y,t,s)\left(\phi(\zeta,z,x,s,u)-\phi(\zeta,y,x,s,t)\right).
$$
Evaluating this at $s=0, x=0$ we get
$$\phi(\zeta,z,y,t,u)=1+\phi(\zeta,0,y,t,0)\left(\phi(\zeta,z,0,0,u)-\phi(\zeta,y,0,0,t)\right).$$
This shows that \rf{BMS} holds for all $n\geq 1$ and $|q|<1$.
Since \rf{BMS} is an identity between the polynomial expressions
in variables $z,y,q$, it must hold for all  $q$.
\end{proof}
We now verify that $\mu_{x,s,t}(dy)$
are the transition probabilities of a Markov process.
\begin{proposition}\label{MMM}
If $0\leq s<t<u$, then $$\mu_{x,s,u}(\cdot)=\int\mu_{y,t,u}(\cdot)\mu_{x,s,t}(dy).$$
\end{proposition}
\begin{proof}
Let
$\nu(dz)=\int\mu_{x,s,t}(dy)\mu_{y,t,u}(dz)$. To show that
$\nu(dz)=\mu_{x,s,u}(dz)$, we verify that $Q_n(z|x,s,u)$ are
orthogonal with respect to $\nu(dz)$. Since $Q_n(z|x,s,u)$ satisfy
the three-step recurrence \rf{Q-rec-G}, we need only to show that
for $n\geq1$ these polynomials integrate to zero. Since $\int Q_k(z|y,t,u)\mu_{y,t,u}(dz)=0$ for $k\geq 1$, by \rf{I} we
have
$$\int Q_n(z|x,s,u) \nu(dz)$$
$$=
\sum_{k=0}^n \left[ \begin{array}{c}n \\ k \end{array}\right]_q \int \left(\int Q_k(z|y,t,u)\mu_{y,t,u}(dz)\right)Q_{n-k}(y|x,s,t)\mu_{x,s,t}(dy)
$$
$$
=\int  Q_{n}(y|x,s,t)\mu_{x,s,t}(dy)=0,
$$
as $n\geq 1$.
\end{proof}
Let  $(X_t)$ be a Markov process with the transition
probabilities  defined for $0\leq s <t$ by
\begin{equation}\label{M-trans}
P_{s,t}(x,dy)=\mu_{x,s,t}(dy),
\end{equation}
where  $\mu_{x,s,t}(dy)$ is the distribution
 orthogonalizing the polynomials  \rf{Q-rec-G}, $X_0=0$.
Since the distribution of $X_t$ is $\mu_{0,0,t}(dx)$,
 the monic polynomials $p_n(x,t)$  orthogonal with
respect to the distribution of $X_t$ are $p_n(x,t)=Q_n(x|0,0,t)$.
These polynomials satisfy a somewhat simpler three-step recurrence
\begin{equation}\label{rec-G}
xp_n(x,t)=p_{n+1}(x,t)+\theta [n]_qp_n(x,t)+(t+\tau[n-1]_q)[n]_q p_{n-1}(x,t), \; n\geq 1.
\end{equation}
Identity \rf{BMS}  can be re-written as
\begin{equation}\label{BMS+}
Q_n(y|x,s,t)=\sum_{k=1}^n
B_{n-k}(x)\left(p_k(y,t)-p_k(x,s)\right),
\end{equation}
where $B_{k}(x)$ are polynomials in variable $x$ such that
$B_0=1$.

If $-1\leq q <1$, then the coefficients of the recurrence
\rf{rec-G} are uniformly bounded. Therefore,  the distribution of $X_t$ has
bounded support, see \cite[Theorem 69.1]{Wall}. If $q=1$, these
are classical Meixner polynomials (see \cite[Ch. VI.3]{Chihara} or
\cite[Sections 4.2 and 4.3]{Schoutens00}), and their distributions
have analytic characteristic functions. This implies that
polynomials are dense in $L_2(X_s,X_u)$, see \cite[Theorem
3.1.18]{Dunkl-Xu}.

We use these observations to extend \cite[(4.4)]{Schoutens00} to
some non-L\'evy processes.
\begin{proposition}\label{LM2} If 
$(X_t)$ is the Markov process with transition probabilities \rf{M-trans}
and $X_0=0$,
then  for $t>s$ and $n\geq 0$ we have
\begin{equation}\label{proj}
E(p_{n}(X_t,t)|\calF_{\leq s})=p_{n}(X_s,s).
\end{equation}
\end{proposition}
\begin{proof} Notice that for $n\geq 1$ we have
$E(Q_n(X_t|X_s,s,t)|X_s)=0$, as $Q_n(y|x,s,t)$ is orthogonal to $Q_0=1$ under the conditional probability \rf{M-trans}.
We use this to prove \rf{proj} by induction.


Since $p_0=1$, \rf{proj} holds true for $n=0$.
Suppose \rf{proj} holds true for all $0\leq n\leq N$. From \rf{BMS+} and the induction assumption
it follows that
$$0=E(Q_{N+1}(X_t|X_s,s,t)|X_s)=B_0(X_s)\left(E(p_{N+1}(X_t,t)|X_s)-p_{N+1}(X_s,s)\right).$$
Since $B_0=1$, this proves that $E(p_{N+1}(X_t,t)|X_s)=p_{N+1}(X_s,s)$, which by the Markov property
implies \rf{proj} for $n=N+1$.
\end{proof}

\begin{proposition}\label{LM4} If $-1\leq q\leq 1$ and $(X_t)$ is the Markov process with transition probabilities \rf{M-trans}
and $X_0=0$,
then \rf{cov}, \rf{LR},
\rf{CV},
and  \rf{q-Var} hold true.
\end{proposition}
\begin{proof} Let $p_n(x,t)$ be the monic polynomials which are
 orthogonal  with respect to the distribution of $X_t$.
For the first part of the proof we will write
 their three step recurrence \rf{rec-G}
as
\begin{equation}\label{rec-0}
xp_n(x,t)=p_{n+1}(x,t)+a_n(t)p_n(x,t)+b_n(t)p_{n-1}(x,t),
\end{equation}
where the coefficients are
\begin{equation}\label{a_n}
a_n(t)=\theta [n]_q, \: b_n(t)=(t+\tau[n-1]_q)[n]_q.
\end{equation}
We will also use the notation
\begin{equation}\label{linear}
a_n(t)=\alpha_n+t \beta_n,\, b_n(t)=\gamma_n+t\delta_n.
\end{equation}

Recall that
\begin{equation}\label{@}
E(p_{n+1}^2(X_t,t))=b_{n+1}(t)E(p_{n}^2(X_t,t)),
\end{equation}
see  \cite[page 19]{Chihara}.

We first verify \rf{CV}. Since $p_1(x,t)=x, p_2(x,t)=x^2-\theta x -t$, from \rf{proj}
we get $E(X_t^2|X_s)=E(p_2(X_t,t)|X_s)+\theta E(p_1(X_t,t)|X_s)+t=
p_2(X_s,s)+\theta p_1(X_s,s)+t=X_s^2+t-s$.

Condition \rf{cov} holds true as $E(X_t)=E(p_1(X_t,t)p_0(X_t,t))=0$, and for $s<t$ we have 
$E(X_sX_t)=E(X_sp_1(X_s,s))=E(p_2(X_s,s)+\theta p_1(X_s,s)+s)=s$.

To verify \rf{LR}, we use the fact that 
 polynomials are dense in $L_2(X_s,X_u)$.
Thus by the Markov property to prove \rf{LR} we
only need to verify
that
\begin{eqnarray}\label{LR-LR}
&E\left(p_n(X_s,s)X_tp_m(X_u,u)\right)&\\ \nonumber
=&a E\left(X_sp_n(X_s,s)p_m(X_u,u)\right)+
bE\left(p_n(X_s,s)X_up_m(X_u,u)\right)&
\end{eqnarray}
for all $m,n\in\sN$ and $0<s<t$. To prove this, we invoke
Proposition \ref{LM2}.
By \rf{proj}
$$
E(p_n(X_s,s) X_t p_m(X_u,u))=E(p_n(X_s,s)X_tp_m(X_t,t)).
$$
Then by  \rf{rec-0} and again using  \rf{proj} we get that the
left hand side of \rf{LR-LR} is
$$
E(p_n(X_s,s)p_{m+1}(X_s,s))+a_m(t)E(p_n(X_s,s)p_m(X_s,s))$$
$$+b_m(t)E(p_n(X_s,s)p_{m-1}(X_s,s)).
$$
Thus the left hand side of the equation is zero, except when
$n=m+1$, $n=m$, or $n=m-1$.

Similar argument applies to the right hand side of \rf{LR-LR}. Thus,
writing  $Ep_m^2$ for $E(p_m^2(X_s,s))$,  equation \rf{LR-LR} takes the form $0=0$,
except for the following three cases.
\begin{enumerate}
\item  Case $n=m+1$. Then the equation reads
$$Ep_{m+1}^2=a(s,t,u) b_{m+1}(s)Ep_m^2+b(s,t,u) Ep_{m+1}^2.$$
By \rf{@}
this holds true as $a+b=1$, see
\rf{a+b}.
\item Case $n=m$.  Then the equation reads
$$a_m(t)Ep_m^2=a(s,t,u) a_m(s)Ep_m^2+b(s,t,u) a_m(u)Ep_m^2.$$
By \rf{a+b}, this equation holds true for any three step recurrence \rf{rec-0} with the coefficients $a_n(t)$ that are
linear in variable $t$.
\item Case $n=m-1$.  In this case, \rf{LR-LR} reads
$$
b_m(t)Ep_{m-1}^2=a(s,t,u)Ep_m^2+b(s,t,u) b_m(u)Ep_{m-1}^2.$$
By \rf{@} this is equivalent to $b_m(t)=a(s,t,u) b_m(s)+b(s,t,u) b_m(u)Ep_{m-1}^2$, which by \rf{a+b} holds true for any
three step recurrence \rf{rec-0} with the coefficients $b_n(t)$ that are linear in variable $t$.
\end{enumerate}

The proof of \rf{q-Var} follows the same plan.
We verify that \rf{QV} holds true with
the parameters given by formulas \rf{A'}, \rf{B'}, \rf{C'}, \rf{D'}, \rf{a'}, \rf{b'}.
(In fact, our proof indicates also how these formulas and the
 recurrence \rf{rec-G} were initially derived.)
To do so, from the three step recurrence \rf{rec-0} we derive
\begin{eqnarray} \label{five-step}
&x^2p_{n-1}(x)=p_{n+1}(x)+(a_n+a_{n-1})p_n(x)& \\ \nonumber +
&(a_{n-1}^2+b_n+b_{n-1})p_{n-1}(x) +
 b_{n-1}(a_{n-1}+a_{n-2})p_{n-2}(x)+b_{n-1}b_{n-2}p_{n-3}(x)&
\end{eqnarray}
 for $n\geq 2$. (Recall that we use the convention $p_{-1}(x)=0$.)

We need to prove that for any $n,m\in\sN$ and $0<s<t$
\begin{eqnarray} \label{QV-QV}
&E\left(p_n(X_s,s)X_t^2p_m(X_u,u)\right)& \\
\nonumber
= &A E\left(X_s^2p_n(X_s,s)p_m(X_u,u)\right) +B
E\left(X_sp_n(X_s,s)X_up_m(X_u,u)\right)& \\ \nonumber
&+C E\left(p_n(X_s,s)X_u^2p_m(X_u,u)\right)
+ \alpha E\left(X_sp_n(X_s,s)p_m(X_u,u)\right)&\\ \nonumber
&+
\beta E\left(p_n(X_s,s)X_up_m(X_u,u)\right) +DE\left(p_n(X_s,s)p_m(X_u,u)\right)
.&\end{eqnarray}
For the remainder of the proof, all the polynomials are evaluated at $(X_s,s)$.
 Using \rf{five-step}, \rf{rec-0} and \rf{proj}, we get
$$
Ep_np_{m+2}+(a_{m+1}(t)+a_m(t))Ep_np_{m+1}+
(a_m^2(t)+b_{m+1}(t)+b_m(t))Ep_np_m$$
$$+
b_m(t)(a_m(t)+a_{m-1}(t))Ep_np_{m-1}+b_m(t)b_{m-1}(t)Ep_np_{m-2}$$
$$=A(Ep_{n+2}p_m+(a_{n+1}(s)+a_n(s))Ep_{n+1}p_m+(a_n^2(s)+b_{n+1}(s)+b_n(s))Ep_np_m
$$
$$+b_n(s)(a_n(s)+a_{n-1}(s))Ep_{n-1}p_m+b_n(s)b_{n-1}(s)Ep_{n-2}p_m
)
$$
$$+
B E\left((p_{n+1}+a_n(s)p_n+b_n(s)p_{n-1})(p_{m+1}+a_m(u)p_m+b_m(u)p_{m-1})\right)
$$
$$+
C(Ep_{n}p_{m+2}+(a_{m+1}(u)+a_m(u))Ep_{n}p_{m+1}+(a_m^2(u)+b_{m+1}(u)+b_m(u))Ep_np_m
$$
$$+b_m(u)(a_m(u)+a_{m-1}(u))Ep_{n}p_{m-1}$$
$$+b_m(u)b_{m-1}(u)Ep_{n}p_{m-2}
)
+ \alpha(Ep_{n+1}p_m+a_n(s)Ep_np_m+b_n(s)Ep_{n-1}p_m)$$$$+
\beta(Ep_np_{m+1}+a_m(u)Ep_np_m+b_m(u)Ep_np_{m-1})
+D Ep_np_m.
$$
Thus the equation \rf{QV-QV} takes the form $0=0$, except for the
following five cases:
\begin{enumerate}
\item Case $n=m+2$.  In this case, equation \rf{QV-QV} reads
$$ E p_{m+2}^2=A b_{m+2}(s)b_{m+1}(s)E p_m^2+Bb_{m+2}(s)Ep_{m+1}^2+CEp_{m+2}^2.$$
By \rf{@}, this is equivalent to \rf{1}, which holds true by our choice of $A,B,C$.
\item Case $n=m+1$.
In this case, equation \rf{QV-QV} reads
$$
\left(a_{m+1}(t)+a_m(t)\right)Ep_{m+1}^2$$
$$=
A b_{m+1}(s)(a_{m+1}(s)+a_m(s))Ep_m^2+B(a_{m+1}(s) Ep_{m+1}^2+b_{m+1}(s)a_m(u)Ep_m^2)$$
$$
+C(a_{m+1}(u)+a_m(u))Ep_{m+1}^2+\alpha b_{m+1}(s)Ep_m^2+\beta Ep_{m+1}^2.
$$
By  \rf{linear} and \rf{@}, this reduces to equation $(\beta_n+\beta_{n-1})=\frac{(u-s)B}{D} \beta_{n-1}$,
 which holds true since $\beta_n=0$, see \rf{a_n}.
\item Case $n=m$.
In this case, equation \rf{QV-QV} reads
$$ \left(a_m^2(t)+b_{m+1}(t)+b_m(t)\right)Ep_m^2=
A \left(a_m^2(s)+b_{m+1}(s)+b_m(s)\right)Ep_m^2+
$$
$$B\left(Ep_{m+1}^2 +a_m(s)a_m(u)Ep_m^2+b_{m}(s)b_m(u)Ep_{m-1}^2\right)
$$
$$+C\left(a_m^2(u)+b_{m+1}(u)+b_m(u)\right)Ep_m^2
+\alpha a_m(s)Ep_m^2+\beta a_m(u) Ep_m^2+DEp_m^2.
$$
After a calculation, this reduces to equation $\delta_{n}+\delta_{n-1}=\delta_{n-1}\frac{(u-s)B}{D}+1$.
The latter
holds true  by \rf{a_n} and \rf{D}.
\item Case $n=m-1$.
In this case, equation \rf{QV-QV} reads
$$
b_m(t)(a_m(t)+a_{m-1}(t))Ep_{m-1}^2=
Ab_m(s)(a_m(s)+a_{m-1}(s))Ep_{m-1}^2$$
$$+
B(a_m(u)Ep_m^2+a_{m-1}(s)b_{m}(u)Ep_{m-1}^2)+
Cb_m(u)(a_m(u)+a_{m-1}(u))Ep_{m-1}^2$$$$+\alpha Ep_m^2+\beta b_m(u) Ep_{m-1}^2.
$$
After a calculation, this reduces to equation
$$(\alpha_{n-1}+\alpha_{n-2})\delta_{n-1}=(1+q)\delta_{n-1}\alpha_{n-2}+
\delta_{n-1}\frac{s\alpha +u \beta }{D(s,t,u)}.$$
The latter
holds true for all $n\geq 2$ by \rf{a_n} and \rf{3}.
\item Case $n=m-2$.
In this case, equation \rf{QV-QV} reads
$$
b_m(t)b_{m-1}(t)Ep_{m-2}^2=A Ep_{m}^2 + B b_m(u)Ep_{m-1}^2+C b_m(u)b_{m-1}(u)Ep_{m-2}^2.
$$
After a calculation, this reduces to equation
$$
\delta_{n-1}\gamma_{n-2}+\delta_{n-2}\gamma_{n-1}=(1+q)\delta_{n-1}\gamma_{n-2}+\delta_{n-1}\delta_{n-2}\frac{As^2+Bsu+Cu^2-t^2}{D}
.$$
Using relation \rf{tau}, this gives
\begin{equation*}\label{****}
\delta_{n-2}\gamma_{n-1}=\tau \delta_{n-1}\delta_{n-2}+q\delta_{n-1}\gamma_{n-2}.
\end{equation*}
The latter is satisfied with the initial condition $\gamma_1=0$ whenever
$$
\gamma_n=\tau [n-1]_q \delta_n.
$$
\end{enumerate}
\end{proof}

From Proposition \ref{iii'} we see that the conditional variance of
a stochastic process $(X_t)$ that satisfies \rf{cov},
\rf{LR}, \rf{QV} with $D\ne 0$, \rf{CV}, and which has at least 3-point support
is given by \rf{Var-q} with parameters $-\infty<\theta<\infty,-1<q\leq 1,\tau\geq 0$.

Let $(Y_t)$ be the Markov process with the transition probabilities
 \rf{M-trans} and the same parameters.  By Proposition \ref{LM4}, this process
satisfies
\rf{cov},  \rf{LR}, \rf{CV}, and \rf{q-Var}.

Since processes $(X_t)$ and $(Y_t)$ satisfy \rf{cov}, \rf{LR}, \rf{CV}, and \rf{q-Var} with the same parameters
$q,\theta, \tau$, and the distribution of $(Y_t)$ is determined uniquely by moments,
therefore by Proposition \ref{L1} the processes  have the same
finite dimensional distributions.
 This establishes our main result.
\begin{theorem}\label{T1} Let $(X_t)_{t\geq0}$ be a
separable square integrable stochastic process which satisfies conditions \rf{cov},
\rf{LR}, \rf{QV}, and \rf{CV}. Suppose that
the coefficient $D$ in \rf{QV} satisfies $D(s,t,u)\ne 0$
for  all 
$0\leq s<t<u$,
and that $1, X_t, X_t^2$ are linearly independent for all $t>0$.
Then there are parameters  $-1< q \leq 1$,  $\theta\in \sR$, and
$\tau\geq 0$  such that $(X_t)$ is a Markov process, with the
transition probabilities \rf{M-trans}, $X_0=0$.

Conversely, for any $-1<q\leq 1, \tau\geq 0, \theta\in\sR$,
the Markov process with transition probabilities \rf{M-trans}
satisfies \rf{cov}, \rf{LR}, \rf{QV}, and \rf{CV}.
\end{theorem}

\begin{remark} If $1,X_u,X_u^2$ are linearly dependent, then the coefficients in
\rf{QV} are not unique; in particular, one can modify
$\beta(s,t,u)$ and $C(s,t,u)$ to get $D(s,t,u)=0$ for all $s<t<u$,
and the assumption $D\ne 0$ makes little sense. However, this  can
sometimes be circumvented, see Theorem \ref{q-Wiener}.
\end{remark}

\begin{remark}\label{R3} 
For $q=1$,  expression \rf{q-Var} depends on the increments of $(X_t)$ only,
i.e., it takes the form
 analyzed in  \cite[Theorem 1]{Wesolowski93}, see also Theorem
\ref{Wesolowski+}.  It is tempting to use this case as a model and
define the $q$-generalizations of the five types of L\'evy
processes determined in  \cite{Wesolowski93}:
\begin{enumerate}
\item
$q$-Wiener processes: $\tau=0$, $\theta=0$.
\item 
$q$-Poisson type processes: $\tau=0$, $\theta\ne 0$.
\item 
$q$-Pascal  type processes: $\tau>0$,  $\theta^2>4\tau$.
\item 
$q$-Gamma type processes: $\tau> 0$, $\theta^2=4\tau$.
\item 
$q$-Meixner type processes: $\theta^2<4\tau$.
\end{enumerate}
Some of these generalizations have already been  studied in the
non-commutative probability; for the $q$-Brownian motion see
\cite{BKS97}, for the $q$-Poisson process  see \cite{Anshelevich03},
\cite{Neu-Speicher94},
\cite{Saitoh-Yoshida00a},  and the references
therein. Anshelevich \cite[Remark 6]{Anshelevich01} states a
recurrence which is equivalent to \rf{Q-rec-G} for $s=0,x=0$; the
latter, written as \rf{rec-G},
 plays the role in our proof of Theorem \ref{T1}.

However, it is also possible that for $|q|<1$ the differences between these processes are
 less pronounced;
when $q=0$,  the transition probabilities in Theorem \ref{Ansh}
share the continuous component and its discrete components also
admit a common interpretation, dispensing with the "cases". The
case of $q=0$ is especially interesting, as it corresponds to
certain free L\'evy processes. As we already pointed out in the
introduction, all free L\'evy non-commutative processes have
classical Markov versions by \cite[Theorem 3.1]{Biane98}.
\end{remark}

\section{Some special cases and examples}\label{SEC 4}

As we already mentioned in the introduction,
 some of the examples we encounter are classical versions of the non-commutative processes that already have been studied.
It might be useful to clarify  terminology.
A non-commutative (real) process
$(\XX_t)_{t\in[0,\infty)}$ is a family of elements of a unital $*$-algebra $\calA$ equipped with a
state (i.e., normalized positive linear functional)  $\Phi:\calA\to\sC$ such that $\XX^*_t=\XX_t$.
A classical version of a non-commutative process $(\XX_t)$ is
a stochastic process $(X_t)$ such that  for every finite
choice $0\leq t_1\leq t_2\leq \dots\leq t_k$ the corresponding moments match:
\begin{equation}\label{EQ: EC1}
\Phi(\XX_{t_1}\dots \XX_{t_k})=E(X_{t_1}\dots X_{t_k}).
\end{equation}
If $\sum a^n\Phi(\XX_{t}^{2n})/{2n}!<\infty$ for some $a>0$, i.e.,  $X_t$ has finite exponential moments, this condition determines uniquely
the finite-dimensional distributions of $(X_t)$. Of course, the left hand side of \rf{EQ: EC1} depends on the order of $\{t_j\}$, which cannot be permuted.

\subsection{$q$-Brownian process}
%
%

For $-1\leq q\leq 1$, the  classical version of
the $q$-Brownian motion, see  \cite[Definition 3.5 and Theorem 4.6]{BKS97}, is a  Markov process
with
the transition probabilities $P_{s,t}(x,dy)$ for $0< s<t$  given  by
\begin{equation}\label{46}
\begin{cases}\displaystyle
 \frac12\left(1+\sqrt{s/t}\right)\delta_{x\sqrt{t/s}}\,(dy)+ \frac12\left(1-\sqrt{s/t}\right)\delta_{-x\sqrt{t/s}}\,(dy)& \mbox{ if } q=-1,
\\
\\
\frac{\sqrt{1-q}}{2\pi \sqrt{4t-(1-q)y^{2}}} \prod_{k=0}^{\infty }\frac{
(t-sq^{k})\left( 1-q^{k+1}\right) \left(
t(1+q^{k})^{2}-(1-q)y^{2}q^{k}\right) }{(t-sq^{2k})^{2}-(1-q)
q^{k}(t+sq^{2k})xy +(1-q)(sy^{2}+tx^{2})q^{2k}}dy & \mbox{ if } -1<q<1, \\
\\
\frac{1}{\sqrt{2\pi(t-s)}}\exp\left(-\frac{(y-x)^2}{2(t-s)}\right)\, dy & \mbox{ if } q=1.\\
\end{cases}
\end{equation}
The support consists of two-point $\pm\frac{\sqrt{t}}{\sqrt{s}}x$ when $q=-1$,
and is bounded $|y|<2\sqrt{t}/\sqrt{1-q}$ when $ -1<q<1$.

The univariate distribution of $X_t,t>0$ is given by the transitions $P_{0,t}(0,dy)$ from $X_0=0$,
which are given by
\begin{equation}\label{45}
\begin{cases}
\frac12\delta_{\sqrt{t}}(dy)+\frac12\delta_{-\sqrt{t}}(dy) & \mbox{ if } q=-1, \\ \\
\frac{\sqrt{1-q}}{2\pi \sqrt{4t-(1-q)y^{2}}}\prod_{k=0}^{\infty
}\left( (1+q^{k})^{2}-(1-q)\frac{y^{2}}{t}q^{k}\right) \prod_{k=0}^{\infty
}(1-q^{k+1})& \mbox{ if } -1<q<1, \\ \\
\frac{1}{\sqrt{2\pi t}}\exp(-\frac{y^2}{2t})\, dy & \mbox{ if } q=1.\\
\end{cases}
\end{equation}

The following shows that the $q$-Brownian motion is characterized by
the assumption that conditional variances are  quadratic, coupled with the additional assumption
that for $t<u$
 the conditional variances $\Var(X_t|\calF_{ \geq u})$ are non-random.

\begin{theorem}\label{q-Wiener} Suppose that $(X_t)_{t\geq 0}$ is a
square-integrable separable process  such that
\rf{cov}, \rf{LR}, \rf{QV}, \rf{CV} hold true,
and  in addition
\begin{equation}\label{RV1}
\Var(X_t|\calF_{ \geq u})=\frac{t}{u}(u-t),
\end{equation}
for all $t<u$.
Then there exists $q\in[-1,1]$ such that
\begin{equation}\label{Var-q}
\Var(X_t|\calF_{\leq s}\vee\calF_{\geq u})=\frac{(t-s)(u-t)}{u-qs}\left(\frac{( 1 - q )
      }{
    {( u - s) }^2}(X_u -X_s  )
    (  sX_u  - uX_s)+1\right).
\end{equation}
Moreover, then $(X_t)$ is  Markov with
 transition probabilities  \rf{46} and \rf{45}.

Conversely, a Markov process, $X_0=0$, with
the transition probabilities given by \rf{46} satisfies conditions \rf{LR}, \rf{QV}, \rf{CV}, and \rf{RV1}.
\end{theorem}

\begin{proof}
Formulas \rf{CV} and  \rf{RV-**}  hold true with $\tau=\theta=0$ by
assumption.  The proof of \rf{tmp1} relies only on these two formulas. Therefore,
$E(X_t^4)=(2+q)t^2$ for some $-1\leq q\leq 1$. In particular, $q=-1$ iff  $(E(X_t^2))^2= E(X_t)^4$, i.e.,
$X_t=\pm \sqrt{t}$ with equal
probabilities. We need to consider separately cases $q=-1$ and $q>-1$.

If $q=-1$, the joint moments are uniquely determined from \rf{LR-}.
Namely,  if $s<t$ and
$m$ is odd then $E(X_t^m|\calF_{\leq s})=t^{(m-1)/2} E(X_t|\calF_{\leq s})=t^{(m-1)/2}X_s$. This determines
all mixed moments  uniquely:
if $n_1,\dots,n_k$ are even numbers, $m_1,m_2,\dots,m_\ell$ are odd numbers, $s_1<s_2<\dots<s_\ell$,
 and $\ell$ is even then we have
$$\E \left(X_{t_1}^{n_1}X_{t_2}^{n_2}\dots X_{t_k}^{n_k} X_{s_1}^{m_1}X_{s_2}^{m_2}\dots X_{s_\ell}^{m_\ell}\right)
$$
$$=
\prod_{j=1}^k t_j^{n_j/2}\prod_{j=1}^{\ell/2} \left(s_{2j-1}^{(m_{2j-1}+1)/2}s_{2j}^{(m_{2j}-1)/2} \right).
$$
If  $\ell$ is odd, then $\E \left(X_{t_1}^{n_1}X_{t_2}^{n_2}\dots X_{t_k}^{n_k} X_{s_1}^{m_1}X_{s_2}^{m_2}\dots X_{s_\ell}^{m_\ell}\right)=0$. Since the same holds true for the two-valued
Markov chain, and its conditional variance can be written as \rf{Var-q}, this ends the proof when $q=-1$.

If $-1<q\leq 1$, then $1,X_t,X_t^2$ are linearly independent for
all $t>0$. To apply Theorem \ref{T1} we need to verify that $D(s,t,u)\ne 0$ for all $s<t<u$.
Suppose $D(s,t,u)=0$ for some $0\leq s<t<u$. Inspecting the proof of Lemma \ref{L001} we see that
equations \rf{CV}, \rf{fut1} (which hold true by assumption) and linear independence imply
\rf{a},\rf{b}, \rf{A}, \rf{B}, and \rf{C} with $D=0$.

We now use these values and the value $E(X_sX_t^2X_u)$ to derive a contradiction. Notice that
\rf{LR-} and \rf{LR+} imply that
$E(X_sX_t^2X_u)=s/t E(X_t^4)=(2+q)st$. On the other hand, since $E(X_t^3)=0$ and $D=0$, from \rf{QV} we get
$$E(X_sX_t^2X_u)=AE(X_s^4)+BE(X_s^2X_u^2)+\frac{s}{u}C E(X_u^4).$$
Since $E(X_s^4)=(2+q)s^2$, and $A,B,C$ are given explicitly,
 a calculation shows that this equation holds true only if $(u-t)(t-s)=0$. Thus
$D(s,t,u)\ne 0$ for all $0\leq s<t<u$.

This shows that the assumptions of
Theorem \ref{T1} are satisfied. Theorem  \ref{T1} shows that $X_t$ is Markov with
uniquely determined transition probabilities. Formulas \rf{46} and
\rf{45} give the distribution which orthogonalizes the corresponding
 Al-Salam--Chihara polynomials, see \cite{Askey-Ismail84MAMS}.
\end{proof}

\subsection{L\'evy processes with quadratic conditional variance}
A special choice of the coefficients in \rf{QV} casts the conditional variance as a quadratic
function of
the increments of the process,
\begin{equation}\label{QV-Levy}
\Var(X_t|\calF_{\leq s}\vee \calF_{\geq u})=C_2(X_u-X_s)^2+C_1(X_u-X_s)+C_0,
\end{equation}
where $C_0=C_0(s,t,u),C_1=C_1(s,t,u),C_2=C_2(s,t,u)$ are deterministic functions of $s<t<u$.

As an application of Theorem \ref{T1},
we give the following version of \cite[Theorem 1]{Wesolowski93}.
\begin{theorem}[Wesolowski]\label{Wesolowski+}
Let $(X_t)_{t\geq 0}$ be a square integrable separable stochastic process
such that the conditions \rf{cov}, \rf{LR}, and \rf{QV-Levy} hold true, and $C_2\ne a b$.
If for every $t>0$ the distribution of $X_t$ has at least 3 point support, then
there are numbers $\theta\in\sR, \tau\geq 0$ such that
 the conditional variance \rf{QV-Levy} is given by
\begin{equation}\label{QV-Levy+}
\Var(X_t|\calF_{\leq s}\vee\calF_{\geq u})=
\frac{(u-t)(t-s)}{u-s+\tau}
\left(\tau\frac{(X_u-X_s)^2}{(u-s)^2}+\theta\frac{X_u-X_s}{u-s}+1\right).
\end{equation}
Moreover, one of the following holds:
\begin{enumerate}
\item $\tau=0,\theta=0$, and  $(X_t)$ is the Wiener processes,
$$E(\exp (iu X_t))=\exp(-t u^2/2).$$
\item  $\tau=0,\theta\not =0$, and $(X_t)$ is a Poisson type processes,
 $$E(\exp (iu X_t))=\exp\left(\frac{t}{\theta^2}(e^{iu\theta }-1)-i \frac{ut}{\theta}\right).$$
\item $\tau>0$ and $\theta^2>4\tau$, and $(X_t)$ is a Pascal
(negative-binomial) type process,
 $$E(\exp (iu X_t))=\left(p  e^{-iu\delta_-}+(1-p) e^{-iu\delta_+}\right)^{-t/\tau},$$
where $\delta_{\pm}=\frac{1}{2}(\theta\pm \sqrt{\theta^2-4\tau})$,
$p=\delta_+/(\delta_+-\delta_-)$.

%
%
%
%

\item  $\tau>0$ and $\theta^2=4\tau$, and $(X_t)$ is a Gamma type process,
 $$E(\exp (iu X_t))=\exp\left(-2 i u t/\theta \right)\left(1-i\frac{u \theta}{2}\right)^{-4t /\theta^2}.$$
\item  $\theta^2<4\tau$, and
$(X_t)$ is a Meixner (hyperbolic-secant) type process,
 $$E(\exp (iu X_t))$$
 $$=\exp\left(i\frac{ u \theta t }{2\tau} \right)\left(\cosh (\frac{\sqrt{4\tau-\theta^2 }u}{2})+
i\frac{\theta}{\sqrt{4\tau-\theta^2 }}
\sinh(\frac{\sqrt{4\tau-\theta^2 }u}{2})  \right)^{-t /\tau}.$$
\end{enumerate}
\end{theorem}

\begin{proof}
We verify that the assumptions of Theorem \ref{T1} are satisfied.

Assumption \rf{QV-Levy} implies that  \rf{QV} holds true with parameters
 $A=C_2+a^2, B=2ab-2C_2, C=C_2+b^2$, and $\alpha+\beta=0$.
Therefore, $A+B+C=1$, which together with \rf{0} implies \rf{C}. Since $C_2\ne ab$ is the same as $C\ne b$,
the latter implies
that $D\ne 0$. Thus we can use Lemma \ref{L000}. From \rf{E1} we get
\rf{CV}.
Theorem \ref{T1} implies that $(X_t)$ is a Markov process with the
 transition probabilities which are identified uniquely from their orthogonal
polynomials, see \cite[Ch VI.3]{Chihara}; see also  \cite[Sections 4.2 and
4.3]{Schoutens00}. In particular, $(X_t)$ has independent and homogeneous
increments, with the distribution of $X_{t+s}-X_s \cong X_t$ as listed in the theorem.

From separability, the usual properties of the trajectories of the
Wiener and Poisson processes follow.
\end{proof}

\subsection{Free L\'evy processes  with quadratic conditional variance}
A special choice of the coefficients in \rf{QV} leads to the following conditional variance
%
\begin{eqnarray}\label{0-Var}
&\Var(X_t|\calF_{\leq s}\vee\calF_{\geq u})=&\\
\nonumber
&ab \left(\frac{(X_u-X_s)(sX_u-uX_s)}{u+\tau}
+\tau\frac{(X_u-X_s)^2}{(u-s)^2}+\theta\frac{X_u-X_s}{u-s}+1\right),
&
\end{eqnarray}
where $a,b$ are the coefficients from \rf{LR}. This formula seems hard to separate by natural assumptions
from the general
expression \rf{q-Var}, but the fact that $q=0$ leads to considerable computational simplifications.
Theorem \ref{T1} in this setting takes the following form, with explicit formulas for the
transition probabilities.
\begin{theorem}\label{Ansh}
Let $(X_t)_{t\geq 0}$ be a square integrable separable stochastic process
such that  the conditions \rf{cov}, \rf{LR}, and \rf{0-Var}
hold true.
If for every $t>0$ the distribution of $X_t$ has at least 3 point support, then
$(X_t)$ is a Markov process with the transition probabilities $P_{s,t}(x,dy)$ given for $0\leq s<t$
by the
Stieltjes-Cauchy transform
\begin{eqnarray}
\label{f-transition}
&\displaystyle\int_{\sR}\frac{1}{z-y}P_{s,t}(x,dy)&  \\  \nonumber =&\displaystyle \frac12
\frac{ (t+s+2\tau)(z - x)+(t-s)\theta  -(t-s)
    \sqrt{ (z-\theta)^2 -
        4 (t+\tau) }}{\tau (z-x)^2+\theta(t-s)(z-x)+tx^2+sz^2-(s+t)x z+(t-s)^2 }.&
\end{eqnarray}
The absolutely continuous part 
 of $P_{s,t}(x,dy)$ is given by the density
$$\frac{1}{2\pi}\frac{(t-s)\sqrt{4(t+\tau)-(y-\theta)^2}}
{\tau(y-x)^2+\theta(t-s)(y-x)+tx^2+sy^2-(s+t)xy+(t-s)^2}, 
$$
supported on $(y-\theta)^2<4(t+\tau)$; the singular part is zero, and the discrete part is zero except for the following cases.
\begin{enumerate}
\item $\tau=0, \theta\ne 0$. Then the discrete part of $P_{s,t}(x,dy)$ is non-zero only for
$x=-s/\theta, 0<s<t<\theta^2$ and is then
$$ \frac{1-t/\theta^2}{1-s/\theta^2}\delta_{-t/\theta}.$$
\item $\tau>0$ and $\theta^2>4\tau$. Then the discrete part of $P_{s,t}(x,dy)$ is
non-zero only if $x=y_*(s)$ and is then
$$
\frac{\left(1-\frac{t}{2\tau}\frac{|\theta|-\sqrt{\theta^2-4\tau}}{\sqrt{\theta^2-4\tau}}\right)^+}
{1-\frac{s}{2\tau}\frac{|\theta|-\sqrt{\theta^2-4\tau}}{\sqrt{\theta^2-4\tau}}}\delta_{y_*(t)},
$$ where
$$y_*(t)=\begin{cases}
-t\frac{\theta-\sqrt{\theta^2-4\tau}}{2\tau} & \mbox{ if } \theta>0\\
-t\frac{\theta+\sqrt{\theta^2-4\tau}}{2\tau} & \mbox{ if } \theta<0
\end{cases}.
$$
\end{enumerate}
\end{theorem}
\begin{proof}
From \rf{0-Var} it follows that $D=ab\ne 0$ and $A+B+C=1$. Since $1,X_t,X_t^2$ are linearly independent by
assumption,   from \rf{E1} we deduce
\rf{CV}.
Thus by Theorem \ref{T1}, $(X_t)$ is a Markov process with the
transition probabilities defined by \rf{Q-rec-G}. It remains to
find the Cauchy-Stieltjes transform of the distribution.

It is well known that the Cauchy-Stieltjes transform
$$
G_{x,s,t}(z)=\int_{\sR}\frac{1}{z-y}P_{s,t}(x,dy)
$$
is given by the continued fraction expansion associated with the orthogonal polynomials,
\cite[page 85]{Chihara}.
The initial polynomials are $$Q_0(y)=1,\; Q_1(y)=y-x,\; Q_2(y)=y^2-(x+\theta) y+\theta x-(t-s).$$
For $n\geq 2$, we have
$$yQ_n(y)=Q_{n+1}(y)+\theta Q_n(y)+(t+\tau)Q_{n-1}(y),$$
so for $n\geq 2$ this is a constant-coefficients recurrence.
Thus the corresponding continued fraction
is
$$\displaystyle
G_{x,s,t}(z)=\displaystyle
\frac{1}{\displaystyle z-x-\frac{\displaystyle  t-s}{z-\theta-
\frac{\displaystyle t+\tau}{\displaystyle z-\theta-\frac{t+\tau}{\ddots}}}}
.$$
This gives
$$
G_{x,s,t}(z)=
\frac{1}{\displaystyle z-x-\frac{\displaystyle t-s}{\displaystyle \phi(z)}},
$$
where
$$\phi(z)=\frac{z-\theta+\sqrt{(z-\theta)^2-4(t+\tau)}}{2}$$ solves the quadratic equation
$$
\phi(z)=z-\theta-\frac{t+\tau}{\phi(z)}.$$
The branch of the root should be taken so that the imaginary parts satisfy
$\Im (z) \Im (G_{x,s,t}(z))\leq 0$. This branch should be taken as the regular branch when $\theta>x$
(with the cut from $-\infty$ to $0$), and as the negative of the regular branch when $\theta<x$.

To get the explicit transition probabilities, we use the Stieltjes inversion formula:
$P_{s,t}(x,dy)$ is the weak limit $\lim_{\eps\to0^+} -\frac{1}{\pi}\Im G(y+i\eps)dy$, see
\cite[page 125]{Akhiezer},
\cite[(4.9)]{Chihara},
\cite[(65.4)]{Wall}.
%
%
The calculations are cumbersome but routine, and an equivalent calculation
has been done by several authors, see
\cite[Theorem 2.1]{Saitoh-Yoshida01}, \cite[Theorem 4
]{Anshelevich01}.  To get the answer given above, one relies on Markov property to determine the
 values of $x$ which can be reached from $0$ at time $s$.
\end{proof}

\begin{remark} The transition probabilities from Theorem \ref{Ansh} can be cast into the form
resembling Theorem \ref{Wesolowski+}. Since the continuous part varies smoothly as we vary the parameters,
the main distinctions between the "five" processes are in the presence of the discrete component.
Accordingly, we have the following cases:
\begin{enumerate}
\item $\tau=0,\theta=0$,  and   $(X_t)$ is
 the free Brownian motion with  the law of $X_t$ given by
$$
\frac{1}{2\pi t}\sqrt{4t-x^2}\: 1_{x^2\leq 4t}\,dx,
$$
see \cite[Section 5.3]{Biane98}.
\item   $\tau=0,\theta\not =0$, and $(X_t)$ is a free Poisson type processes
with  the law of $X_t$ given by
$$
\left(1-t/\theta^2\right)^+\delta_{-{t}/{\theta}}(dx)+
\frac{1}{2\pi}\frac{1}{\theta x+t}\sqrt{4t-(x-\theta)^2}\:1_{(x-\theta)^2\leq 4t}\,dx,
$$
compare \cite[Section 2.7]{Voiculescu00}.
\item  $\tau>0$ and $\theta^2>4\tau$, and $(X_t)$ is a free Pascal (Negative binomial) process
with the  law of $X_t$ given by
$$
p_*(t)\delta_{x_*}+
\frac{1}{2\pi}\frac{t}{\tau x^2+t\theta x+t^2}\sqrt{4(t+\tau)-(x-\theta)^2}
\:1_{(x-\theta)^2\leq 4(t+\tau)}\,dx,
$$
where
$$
p_*(t)=
\left(1-\frac{t}{2\tau}\frac{|\theta|-\sqrt{\theta^2-4\tau}}{\sqrt{\theta^2-4\tau}}\right)^+
,
$$
and
$$x_*(t)=\begin{cases}
{{t}(\sqrt{\theta^2-4\tau}-\theta)/(2\tau)} &\mbox{ if } {\theta>0},\\
\\
{-{t}(\sqrt{\theta^2-4\tau}+\theta)/(2\tau)} & \mbox{ if } {\theta<0},
\end{cases}$$
compare \cite[Theorem 4 ]{Anshelevich01}.
\item $\tau>0$ and $\theta^2=4\tau$ and $(X_t)$ is a free Gamma type process
with the  law of $X_t$ given by
$$
\frac{1}{2\pi}\frac{4t}{(x\theta+2t)^2}\sqrt{4t+\theta^2-(x-\theta)^2}\:1_{(x-\theta)^2\leq 4t+\theta^2}\,dx,
$$
compare \cite[Theorem 4 ]{Anshelevich01}.
\item    $\theta^2<4\tau$, and
$(X_t)$ is a free Meixner (hyperbolic-secant) type process
with the law of $X_t$ given by
$$\frac{1}{2\pi}\frac{t}{\tau x^2+t\theta x+t^2}\sqrt{4(t+\tau)-(x-\theta)^2}\:1_{(x-\theta)^2\leq 4t}\,dx,$$
compare the ``Continuous Binomial process" in \cite[Theorem 4 ]{Anshelevich01}.
\end{enumerate}
\end{remark}

We remark that these measures are closed under the free convolution,
 $\mathcal{L}(X_{t+s})=\mathcal{L}(X_{t})\boxplus \mathcal{L}(X_{t})$; this is well known,
 and can be easily seen from the corresponding $R$-series which for $\tau>0$ is
$$
R_{X_t}(z)=t \frac{
1 - z\theta - \sqrt{(1 - z\theta )^2 - 4z^2\tau } }{2z\tau },
$$
compare \cite[Proposition 3.4]{Bozejko-Bryc-04}.
The free
Brownian and free Poisson processes have been studied in considerable detail,
see \cite{Voiculescu00} and the references therein.
Symmetric free Meixner distribution appears in \cite[Theorem 3]{Bozejko-Speicher91}, and
in \cite{Bozejko-Leinert-Speicher}. According to \cite[Theorem 3.2(2)]{Saitoh-Yoshida01}, these
laws are infinitely divisible with
respect to the free convolution, with explicit L\'evy representations.
 All five distributions occur in
Anshelevich \cite[Theorem 4]{Anshelevich01}; Anshelevich also
points out that the correspondence between the classical and free
Levy processes based on the values of parameters $\theta, \tau$
does not match the Bercovici-Pata bijection.

\subsection{Binomial Example}
The coefficients in \rf{LR} and \rf{QV} alone do not determine the distribution of a process, and
\rf{LR} and \rf{QV} may be satisfied
by processes with univariate distributions different than those listed in Theorem \ref{T1}.

\begin{proposition}\label{P-MP}
 Let
$p:[0,\infty)\to[0,\infty)$ be such that
$\int_0^{\infty}p(x)\:dx< 1$. Fix $m\in\sN$ and let $\pi(s,t)=\int_s^t\:p(x)\:dx$. The Markov process $(Y_s)_{s\ge 0}$
with $Y_0=0$ and the transition probabilities
$$
P(Y_t=j|Y_s=i)=\frac{(m-i)!}{(j-i)!(m-j)!}
\frac{\left(\pi(s,t)\right)^{j-i}\left(1-\pi(0,t)\right)^{m-j}}
{\left(1-\pi(0,s)\right)^{m-i}}\;,
$$
for $0\le i\le j\le m$ and any $0\le s<t$,
satisfies \rf{LR} and \rf{QV} with the coefficients that do not depend on the parameter  $m\in\sN$. Namely,
\begin{equation}\label{NB-LR}
E(Y_t|\calF_{\leq s}\vee \calF_{\geq u}
)=\frac{\pi(t,u)}{\pi(s,u)}Y_s+
\frac{\pi(s,t)}{\pi(s,u)}Y_u
\end{equation}
and
\begin{equation}\label{NB-QV}
\Var(Y_t|\calF_{\leq s}\vee \calF_{\geq u})=
\frac{\pi(s,t)\pi(t,u)}{\left(\pi(s,u)\right)^2}(Y_u-Y_s)
.\end{equation}
\end{proposition}

\begin{proof} We first show that the
transition probabilities are consistent. For any $0\le s<t<u$ and
integers $i,n\ge 0, i+n\leq m$
$$
P(Y_u=i+n|Y_s=i)=\sum_{j=0}^n\:P(Y_u=i+n|Y_t=i+j)P(Y_t=i+j|Y_s=i)
$$
$$
=\sum_{j=0}^n \frac{(m-i)!\left(\pi(t,u\:)\right)^{n-j}
\left(1-\pi(0,u\:)\right)^{m-i-n}\left(\pi(s,t\:)\right)^j}
{j!(n-j)!(m-i-n)!\left(1-\pi(0,s\:)\right)^{m-i}}
$$
$$
=
\frac{(m-i)!\left(1-\pi(0,u )\right)^{m-i-n}}
{n!(m-i-n)!\left(1-\pi(0,s )\right)^{m-i}}\sum_{j=0}^n
 \left(\begin{array}{c}n\\j\end{array}\right)
\left(\pi(s,t)\right)^j
\left(\pi(t,u )\right)^{n-j}
$$
$$
=\frac{(m-i)!}{n!(m-i-n)!}\frac{\left(1-\pi(0,u )\right)^{m-i-n}}
{\left(1-\pi(0,s )\right)^{m-i}}\left(\pi(s,t )+\pi(t,u )\right)^n.
$$
$$
=\frac{(m-i)!}{n!(m-i-n)!}\frac{\left(1-\pi(0,u )\right)^{m-i-n}}
{\left(1-\pi(0,s )\right)^{m-i}}\left(\pi(s,u )\right)^n.
$$

Then the joint distribution of $(Y_s,Y_t,Y_u)$ is given by
$$
P(Y_u=i+j+k,\:Y_t=i+j,\:Y_s=i)$$
$$=P(Y_u=i+j+k|Y_t=i+j)P(Y_t=i+j|Y_s=i)P(Y_s=i|Y_0=0)
$$
$$=
 \left(\begin{array}{c}m-i-j\\k\end{array}\right)
\frac{\left(\pi(t,u )\right)^k
\left(1-\pi(0,u )\right)^{m-i-j-k}}{\left(1-\pi(0,t )\right)^{m-i-j}}
$$
$$\times
 \left(\begin{array}{c}m-i\\j\end{array}\right)
\frac{\left(\pi(s,t)\right)^j\left(1-\pi(0,t)\right)^{m-i-j}}
{\left(1-\pi(0,s)\right)^{m-i}}
 \left(\begin{array}{c}m\\i\end{array}\right)
\left(\pi(0,s)\right)^i\left(1-\pi(0,s)\right)^{m-i}
$$
$$
=\frac{m!}{i!j!k!(m-i-j-k)!}\left(\pi(0,s)\right)^i\left(\pi(s,t)\right)^j
\left(\pi(t,u)\right)^k\left(1-\pi(0,u)\right)^{m-i-j-k}.
$$
From this, it is easy to see that
conditionally on $Y_s,Y_u$, the increment $Y_t-Y_s$ has the
binomial distribution with $Y_u-Y_s$ trials and the probability of
success ${\pi(s,t )}/{\pi(s,u )}$, i.e.,
$$P(Y_t=k+i|Y_s=i,Y_u=i+n)= \left(\begin{array}{c}n\\k\end{array}\right)
\left(\frac{\pi(s,t )}{\pi(s,u )}\right)^k
\left(\frac{\pi(t,u )}{\pi(s,u )}\right)^{n-k}.$$
Therefore
$$
E(Y_t|Y_s,Y_u)=Y_s+\frac{\pi(s,t )}{\pi(s,u )}(Y_u-Y_s),
$$
and \rf{NB-LR} follows from the Markov property. Similarly, \rf{NB-QV} is a consequence of Markov property
and the formula for
the variance of the binomial distribution.

\end{proof}

\begin{remark} For $s\leq t$ the conditional distribution of  $Y_t - Y_s$ given $Y_s$ is binomial
$b(m-Y_s,\pi(s,t)/(1-\pi(0,s))$, which gives
$$Cov (Y_s, Y_t) = m \pi(0,s)\left(1-\pi(0,t)\right).$$
\end{remark}
%
%
%
%
%

\subsection*{Acknowledgement}  Part of the research of WB was conducted while visiting
the Faculty of Mathematics and Information Science of
Warsaw University of Technology.
The authors thank M. Bo\.zejko for bringing to their attention
several references, to Hiroaki Yoshida for information pertinent to Theorem \ref{Ansh},
 and to M. Anshelevich,
W. Matysiak, R. Speicher, P. Szab{\l}owski, and  M. Yor for helpful comments and discussions.
Referee's comments lead to several improvements in the paper.

\bibliographystyle{acm} 
\bibliography{q-reg,Vita,Wesol}

\end{document}